\documentclass[12pt]{amsart}

\usepackage{amsfonts,amssymb,amsmath,amsthm,graphicx}
\usepackage[mathscr]{eucal}

\newtheorem{thm}{Theorem}[section]
\newtheorem{lem}[thm]{Lemma}
\newtheorem{conj}[thm]{Conjecture}

\theoremstyle{remark}

   \def\Int{\operatorname{Int}} 
   \def\pr{\operatorname{pr}} 
   \def\tb{\operatorname{tb}}   
   \def\TB{\operatorname{TB}}   
   \def\a{\alpha} \def\d{\delta}
   \def\e{\varepsilon}  
    
   \def\phi{\varphi}  \def\s{\sigma} \def\t{\tau} 
   
   \def\G{\Gamma}
     \def\Bd{\partial} 
     \def\Nb#1#2{N_#1(#2)}
     \def\C{\mathbb C} \def\F{\mathbb F} \def\N{\mathbb N} 
     \def\R{\mathbb R} \def\Z{\mathbb Z}
     \def\sub{\subset} \def\emptyset{\varnothing}
     \def\As{\mathscr A}  
      \def\Hs{\mathscr H} 
       \def\Ls{\mathscr L} 
     
     \def\P{\mathbb P} \def\Ps{\mathscr P} 
     \def\Ls{\mathscr L} 
     \def\CLs{\C\thinspace\Ls}
      
      \def\Ts{\mathscr T}

     \def\A#1{\As(#1)}\def\H#1{\Hs(#1)}
     \def\Kinf#1{K_\infty({#1})}
     
\makeatletter
\mathchardef\Unfold@="1007 
 \def\Unfold{\DOTSB\Unfold@\slimits@}
\makeatother
\begin{document}

\title[Links at infinity of hyperbolic line arrangements]%
{Some fibered and non-fibered \\
links at infinity of hyperbolic \\
complex line arrangements}
\author{Lee Rudolph}
\address{Department of Mathematics, Clark University, 
Worcester MA 01610 USA}  
\thanks{Partially supported by NSF (DMS-9504832) and by 
the Fonds National Suisse.}
\email{lrudolph@black.clarku.edu}            

\keywords{
Arrangement,  
complex hyperbolic plane,
divide,
fibered link,
Legendrian knot,
link at infinity,
Murasugi sum,
quasipositivity,
unfolding}

\subjclass{Primary 57M25, 52C30, 51M99}

\begin{abstract}
Let $\F$ be $\R$ or $\C$, $d:=\dim_\R(\F)$.  
Denote by $\Ps(\F)$ either the affine plane $\A\F$ 
or the hyperbolic plane $\H\F$ over $\F$.  
An arrangement $\Ls$ of $k$ lines in $\Ps(\F)$ (pairwise 
non-parallel in the hyperbolic case) has a link at infinity 
$\Kinf\Ls$ comprising $k$ unknotted $(d-1)$-spheres in
$S^{2d-1}$, whose topology reflects the combinatorics 
of $\Ls$ ``at infinity''.  The class of links at infinity 
of affine $\F$-line arrangements is properly included in the 
class of links at infinity of hyperbolic $\F$-line arrangements.
Many links at infinity of (essentially non-affine) connected
hyperbolic $\C$-line arrangements are far from being fibered.  
In contrast, if the (affine or hyperbolic) $\R$-line arrangement 
$\Ls_\R\sub\Ps(\R)$ is connected, and $\Ls=\CLs_\R\sub\Ps(\C)$ 
is its complexification, then $\Kinf\Ls$ is fibered.
\end{abstract}
\maketitle

\section{Introduction; statement of results}\label{intro}

A link $Q\sub S^3$ is \emph{split} if 
$\pi_2(S^3\setminus Q)\ne\{0\}$, 
and \emph{fibered} if $Q$ has an \emph{open book}, i.e., 
a map $f: S^3\to\C$ with $Q=f^{-1}(0)$ such that 
$0$ is a regular value and $f/|f|:S^3\setminus Q\to S^1$ 
is a fibration.  A fibered link is not split.
 
\textbf{Theorem \ref{Theorem A}.}\emph{
For every knot $K\sub S^3$, 
for every $t\in\Z$ not greater than 
the maximal Thurston-Bennequin invariant $\TB(K)$,
for all sufficiently large $k\in\N$, 
there is a connected hyperbolic $\C$-line arrangement 
whose link at infinity is a \hbox{$t$-twisted} \hbox{$+$-clasped}
chain of $k$ unknots of type $K$.
This link is not split, and it is fibered iff $K=O$ and $t=-1$.
}%

\textbf{Theorem \ref{Theorem B}.}\emph{
The link at infinity of the complexification 
$\C\Ls_\R\sub\P(\C)$ of an $\R$-line arrangement 
$\Ls_\R\sub\P(\R)$ is fibered iff it is not split 
iff $\Ls_\R$ is connected.
}%

Section~\ref{hyperbolic line arrangements} 
is an introduction to hyperbolic line arrangements 
and their links at infinity.
Theorem~\ref{Theorem A} is proved in Section~\ref{proof of Theorem A} 
by ``Legendrian inscription''.
Theorem~\ref{Theorem B} is proved in Section~\ref{proof of Theorem B}
as an application of A'Campo's method of \emph{divides} 
\cite{A'Campo:divides97,A'Campo:divides00}.  
Section~\ref{remarks} contains various remarks.

\section{Hyperbolic line arrangements and 
their links at infinity}\label{hyperbolic line arrangements}

\subsection{Affine and hyperbolic lines and planes}

Let $\F$ be $\R$ or $\C$; let $d:=\dim_\R(\F)$.  
Denote the euclidean norm on the $\R$-vectorspace $\F^2$ by $\|\cdot\|$.
Coordinatize the affine plane over $\F$ in the usual way:
that is, let $\A\F:=\F^2$, with its standard structure of 
$\F$-analytic manifold, and let an \emph{affine $\F$-line} 
be any translate $(\s,\t)+V$ by $(\s,\t)\in\F^2$
of a $1$-dimensional $\F$-subvectorspace $V\sub\F^2$.
The \emph{complexification} of 
the affine $\R$-line $L=(\s,\t)+V$ is 
the affine $\C$-line $\C L:=(\s,\t)+\C V$.

Coordinatize the hyperbolic plane over $\F$ by its Klein 
model (convenient here, where the hyperbolic metric 
plays no significant role): that is, let 
$\H\F:=\{(\s,\t)\in \F^2: \|(\s,\t)\|<1\}=\Int D^{2d}$, with 
the $\F$-analytic manifold structure induced from $\F^2$, and
declare that a \emph{hyperbolic $\F$-line} is any 
non-empty intersection of $\H\F$ with an affine $\F$-line.
The unique affine $\F$-line containing the hyperbolic
$\F$-line $L\sub\H\F$ will be called the \emph{extension} 
of $L$, denoted $L^e$.  The \emph{complexification} of
the hyperbolic $\R$-line $L$ is $\C L:=\H\C\cap \C(L^e)$;
of course $(\C L)^e=\C(L^e)$.

As usual in the Klein model, 
distinct hyperbolic $\F$-lines $L_1$ and $L_2$ 
are \emph{concurrent} if $L_1^e\cap L_2^e\cap\Int D^{2d}\ne\emptyset$,
\emph{parallel} if $L_1^e\cap L_2^e\cap\Bd D^{2d}\ne\emptyset$,
and \emph{hyperparallel} if $L_1^e\cap L_2^e\cap D^{2d}=\emptyset$.

\subsection{Affine and hyperbolic line arrangements}
For present purposes, 
an \emph{affine $\F$-line arrangement} 
is the union $\Ls=\bigcup L_i$ of a (finite) family
of pairwise distinct affine $\F$-lines $L_i$.
A \emph{hyperbolic $\F$-line arrangement} 
is the union $\Ls=\bigcup L_i$ of a (finite) family
of pairwise distinct, pairwise non-parallel 
hyperbolic $\F$-lines $L_i$.  In either case, 
a \emph{node} of $\Ls$ of \emph{degree $m\ge 2$} 
is a point contained in $m$ lines $L_i\sub\Ls$; $\Ls$ is
\emph{in general position} if it has no nodes of degree $m>2$.  

Let $\Ls$ be a hyperbolic $\F$-line arrangement.
The \emph{extension} of $\Ls$ is the affine 
$\F$-line arrangement $\Ls^e:=\bigcup L_i^e\subset \A\F$.
Easily, every affine $\F$-line arrangement in $\A\F$ is 
homothetic to $\Ls^e$ for some hyperbolic $\F$-line arrangement 
$\Ls$. In case $\F=\R$, the \emph{complexification} of $\Ls$ is 
the hyperbolic $\C$-line arrangement $\CLs:=\bigcup \C L_i$.
The \emph{boundary} (resp., \emph{closure}) of $\Ls$ is
$\Bd\Ls:=\Ls^e\cap S^{2d-1}$ (resp., $\overline\Ls:=\Ls^e\cap D^{2d}$);
thus $\Bd\Ls$ is the union of a family of $k$ pairwise 
disjoint $(d-1)$-spheres in $S^{2d-1}$.   

In case $\F=\C$, 
$\Bd\Ls$ is endowed with a canonical orientation: the orientation 
of the component $L_i^e\cap S^3\sub\Bd\Ls$ is that induced on it by 
the canonical (complex) orientation of $L_i^e$.
In case $\F=\R$, 
$\Bd\Ls$ is endowed with a canonical partition into 
$2$-sets, namely, $L_i^e\cap S^1$.

\begin{def}\label{link at infinity}
Any submanifold of $S^{2d-1}$ isotopic to $\Bd\Ls$, 
by an ambient isotopy respecting the canonical 
orientation \emph{(}resp., partition\emph{)} in case 
$\F=\C$ \emph{(}resp., $\F=\R$\emph{)}, will be called 
the \emph{link at infinity} of $\Ls$ and denoted $\Kinf\Ls$.  
\end{def}

Call a hyperbolic $\F$-line arrangement $\Ls$ \emph{essentially
affine} if 
$\mbox{$\Ls^e\setminus \Ls$}\subset\mbox{$\F^2\setminus\H\F$}$
is a collar of $\Bd\Ls$ (equivalently, if every node of
$\Ls^e$ is already a node of $\Ls$).  Note that a hyperbolic 
$\R$-line arrangement is essentially affine iff its 
complexification is essentially affine.  Evidently,
affine $\F$-line arrangements are combinatorially and 
topologically precisely the same as essentially affine 
hyperbolic $\F$-line arrangements; in particular, the link 
at infinity of an affine $\F$-line arrangement $\Ls$ is
well-defined (even though $\Bd\Ls$ does not exist), and 
the notation $\Kinf\Ls$ will be used in this case also.

The ambient isotopy type of $\Kinf\Ls$ is readily seen to 
be invariant under a number of operations on $\Ls$, such as 
topological equivalence, isotopy, and homotopy (suitably defined), 
including perturbations which are ``sufficiently small''.  
One such invariance in particular will be useful below; 
its truth depends on the exclusion of parallel lines from 
hyperbolic $\F$-line arrangements.

\begin{lem}\label{perturb}
Let $\Ls$ be a hyperbolic $\F$-line arrangement. 
For almost all arrangements $\Ls'$ which are sufficiently 
close to $\Ls$ in the euclidean metric on $\H\F$, 
$\Ls'$ is in general position and $\Kinf{\Ls'}$ is ambient
isotopic to $\Kinf{\Ls}$.
\qed
\end{lem}

\section{A construction of non-fibered 
hyperbolic $\C$-line arrangements}\label{proof of Theorem A}

\subsection{Legendrian knots and the Thurston-Bennequin invariant
of a knot}\label{about Legendrian knots}

At each point of $S^3\sub\C^2$, the tangent (real) $3$-space
to $S^3$ contains a unique affine $\C$-line, 
the \emph{contact plane} at that point.  
An unoriented knot $K\sub S^3$ whose tangent line
at each point is contained in the contact plane at that point
is called \emph{Legendrian}.
A Legendrian knot $K$ is framed by a canonical normal linefield
(the orthogonal complement of its tangent linefield in the
contact planefield); the integer associated in the usual way 
to that framing is the \emph{Thurston-Bennequin invariant} 
$\tb(K)$.  For any knot $K\sub S^3$, let
$\TB(K):=\sup\{\tb(K'): K'\text{ is Legendrian and ambient 
isotopic to }K\}$ denote the \emph{maximal Thurston-Bennequin 
invariant} of $K$.  It is known that $\TB(K)$ is an integer (i.e., 
$\TB(K)\neq\pm\infty$), and that if $\TB(K)\ge t\in\Z$, then 
there is a Legendrian knot $K'$ ambient isotopic to $K$
with $\tb(K')=t$.  (See \cite{Bennequin} and references
cited therein.)

\begin{figure}
\begin{center}

\includegraphics[width=5in]{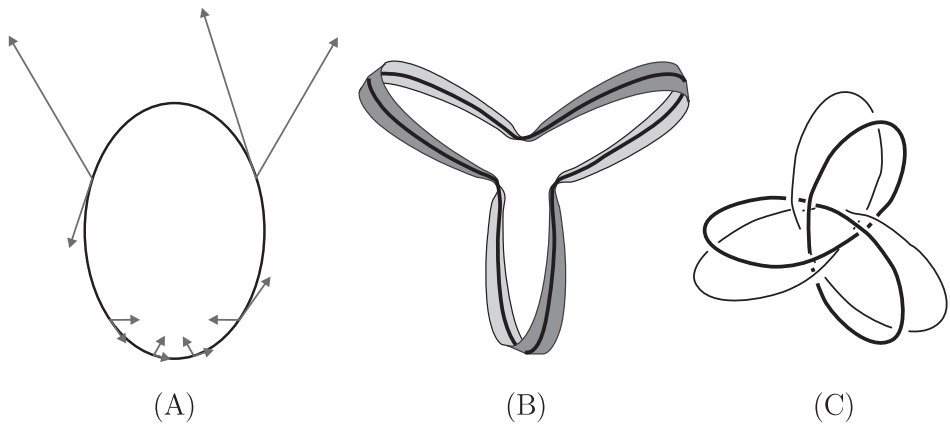}
\end{center}
\caption{(A)~A Legendrian unknot $O\{1,-1\}$ decorated
with some tangent vectors and their corresponding 
Legendrian normal vectors.
(B)~A Legendrian unknot $O\{1,-3\}$ as the core circle 
of an annulus $A(O,-3)$.  (C)~A Legendrian negative 
trefoil $O\{2,-3\}$ and a parallel obtained by pushing 
off along the canonical normal linefield; their linking 
number equals $-6=2(-3)$, which is $\TB(O\{2,-3\})$ 
(\cite[Theorem~8]{constqp4}).}
\label{legendrian knots}
\end{figure}

Fig.~\ref{legendrian knots} illustrates several Legendrian
knots, projected stereographically via
$$
S^3\setminus\{(0,-\mathrm{i})\}\to\C\times\R:%
(z,w)\mapsto(z,\mathrm{Re}(w))/(1+\mathrm{Im}(w))
$$
and thence orthogonally via $\pr_1$ to $\C$, with their 
canonical framings depicted in various ways.  To simplify the 
calculations attendant on the preparation of Fig.~\ref{legendrian 
knots} (and Fig.~\ref{more chains}), each knot illustrated 
has been taken to be of the particularly simple form 
$S^3\cap \{2^{-1/2}\tau^p,2^{-1/2}\tau^{-q}):\tau\in S^1\}$, 
where $p,q\in\N$.

\subsection{Annular plumbing and chains of unknots}
\label{chains of unknots}

Given a knot $K\sub S^3$ and an integer $t$, denote by 
$A(K,t)\sub S^3$ any annulus such that
$K\sub\Int A(K,t)$, $K$ bounds no disk on $A(K,t)$,
and the Seifert matrix of $A(K,t)$ is $[t]$
(that is, when $A(K,t)$ is endowed with an orientation, 
the linking number in $S^3$ of its two boundary components
is $-t$);
for example, letting $O$ denote an unknot, $A(O,-1)$ is
a \emph{positive Hopf annulus}.

A \emph{transverse arc} of $A(K,t)$ is any arc $\tau\sub A(K,t)$
with $\Bd\tau=\tau\cap\Bd A(K,t)$ which intersects $K$ transversely
and at a single point.  Given pairwise disjoint transverse arcs 
$\tau_i\sub A(K,t)$, $i=1,\dots,k$, there are pairwise disjoint
$3$-disks $N_i\sub S^3$ such that $N_i\cap A(K,t)=\Bd N_i\cap A(K,t)$ 
is a regular neighborhood of $\tau_i$ on $A(K,t)$,
and positive Hopf annuli $A(O_i,-1)\sub N_i$ such that
$A(K,t)\cap\Bd N_i = A(O_i,-1)\cap \Bd N_i = A(O_i,-1)\cap A(K,t)$ 
and $K\cap A(O_i,-1)$ is a transverse arc of $A(O_i,-1)$.  The union
$$
F(K,t,k):=A(K,t)\cup\bigcup_{i=1}^k A(O_i,-1)\sub S^3
$$ 
is an orientable surface; it is, in fact, an iterated \emph{plumbing}
(quadrilateral Murasugi sum) of annuli 
\cite{Gabai:Murasugi1,constqp5}.  

As is easily seen, 
$F(K,t,k)$ equipped with either orientation is ambient isotopic 
to $F(K,t,k)$ equipped with the opposite orientation, so 
its boundary $\Bd F(K,t,k)$ has an orientation which is
well-defined up to ambient isotopy.
For $k\ge 2$, $\Bd F(K,t,k)$ is a link of $k$ components;
each component is an unknot, and any $2$-component sublink 
is either a trivial link (if it is split) or isotopic 
to $\Bd A(O,-1)$ (if it is not split).  
Call $\Bd F(K,t,k)$ a 
\emph{\hbox{$t$-twisted} \hbox{$+$-clasped} chain 
of $k$ unknots of type $K$}.

\begin{figure}
\begin{center}
\includegraphics[width=5in]{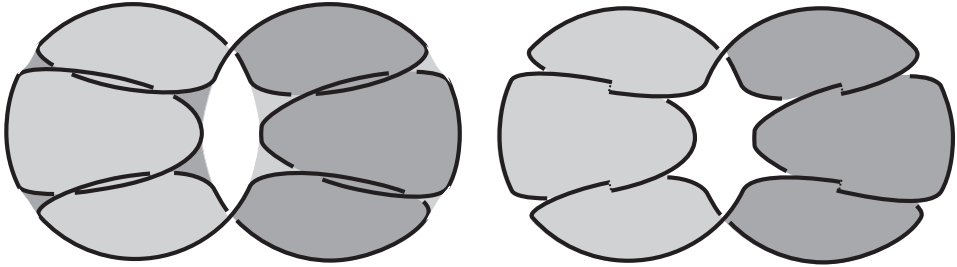}
\end{center}
\caption{Two views of a $(-1)$-twisted $+$-clasped chain 
of $4$ unknots of type $O$.}
\label{chain}
\end{figure}

Fig.~\ref{chain} illustrates two views of 
a \hbox{$(-1)$-twisted} \hbox{$+$-clasped} chain of 
$4$ unknots of type $O$: at the left, it is shown bounding
the iterated plumbing of annuli $F(O,-1,4)$; 
at the right, it is shown bounding the union of four 
individually embedded $2$-disks immersed in $S^3$ 
with (positive) clasp singularities.

\begin{thm}\label{Theorem A}
For every knot $K\sub S^3$, 
for every $t\in\Z$ not greater than $\TB(K)$,
for all sufficiently large $k\in\N$, 
there is a connected hyperbolic $\C$-line arrangement $\Ls$
such that $\Kinf\Ls=\Bd F(K,t,k)$.
This link is fibered iff $K=O$ and $t=\TB(O)=-1$.
\end{thm}
 
\begin{proof} Without loss of generality, $K$ is Legendrian 
and $t=\tb(K)$.
Given $p,q\in S^3$ with $p\ne q$, let $M(p,q)\sub\C^2$ 
denote the affine $\C$-line through $p$ and $q$;
let $M(p,p)$ denote the contact plane at $p$.
Thus $M(p,q)\cap S^3$ is a (round) circle for $p\ne q$, 
and $M(p,p)\cap S^3$ reduces to the singleton $\{p\}$.
Let $d(p,q)$ denote the (euclidean) diameter of $M(p,q)\cap S^3$.
Because $K$ is Legendrian, the limit of $d(p,q)$ as $q\in K$ 
approaches $p\in K$ along $K$ exists and is $d(p,p)=0$.  
It follows that
for all $\e>0$ there exists $\d>0$ such that
if $p,q\in K$ and $\|p-q\|<\d$, then $d(p,q)<\e$.  
Consequently, for all sufficiently large $k$, 
there is a sequence of $k$ points 
$p_1,\dots,p_k, p_{k+1}=p_1$ on $K$ 
(cyclically ordered identically by their indices and
by their position on $K$) and a piecewise-smooth annulus 
$A(K,t)$ with $A(K,t)\supset M(p_i,p_{i+1})\cap S^3$ for all $i$.
Denote by $\widetilde L_i$ the hyperbolic $\C$-line
$M(p_i,p_{i+1})\cap \H\C$. 
The union $\widetilde\Ls:=\bigcup_{i=1}^k \widetilde L_i$ 
(which is, as it were, an ``inscribed hyperbolic $\C$-polygon'' 
of $K$) is not a hyperbolic $\C$-line arrangement 
as defined in this paper, 
since $\widetilde L_i$ and $\widetilde L_{i+1}$ are parallel;
but almost any small perturbation of $\widetilde\Ls$ will be 
a hyperbolic $\C$-line arrangement, 
and for a suitable such perturbation $\Ls$ (easily achieved 
by moving all the affine $\C$-lines $M(p_i,p_{i+1})$ slightly 
closer to the origin) the link at infinity $\Kinf\Ls$ is isotopic 
to $\Bd F(K,t,k)$. 

It remains to be shown that, with $K$ and $t$ as above,
$\Bd F(K,t,k)$ is fibered iff $K=O$ and $t=-1$.  
This can be done by invoking results 
about quasipositivity, Murasugi sums, and fiber surfaces.
The annulus $A(K,t)$ is quasipositive (by \cite{constqp4}),
as is each $A(O_i,-1)$, so (by \cite{constqp5}) the plumbed
surface $F(K,t,k)$ is quasipositive and therefore 
a least-genus surface for its boundary.  
As is well-known, a least-genus Seifert surface bounded 
by a fibered link must be a fiber surface.  
By Gabai \cite{Gabai:Murasugi1}, a Murasugi sum 
of Seifert surfaces is a fiber surface iff all of the 
summands are fiber surfaces.  Since $A(O_i,-1)$ is indeed a 
fiber surface, $\Bd F(K,t,k)$ is fibered iff $A(K,t)$ is a fiber 
surface.  The only annuli which are fiber surfaces are $A(O,-1)$ 
and its mirror image $A(O,1)$.  Since $\TB(O)=-1$ 
(cf.\ \cite{classical}), 
$A(K,t)$ must be $A(O,-1)$. 
\end{proof}

\begin{figure}
\begin{center}
\includegraphics[width=5in]{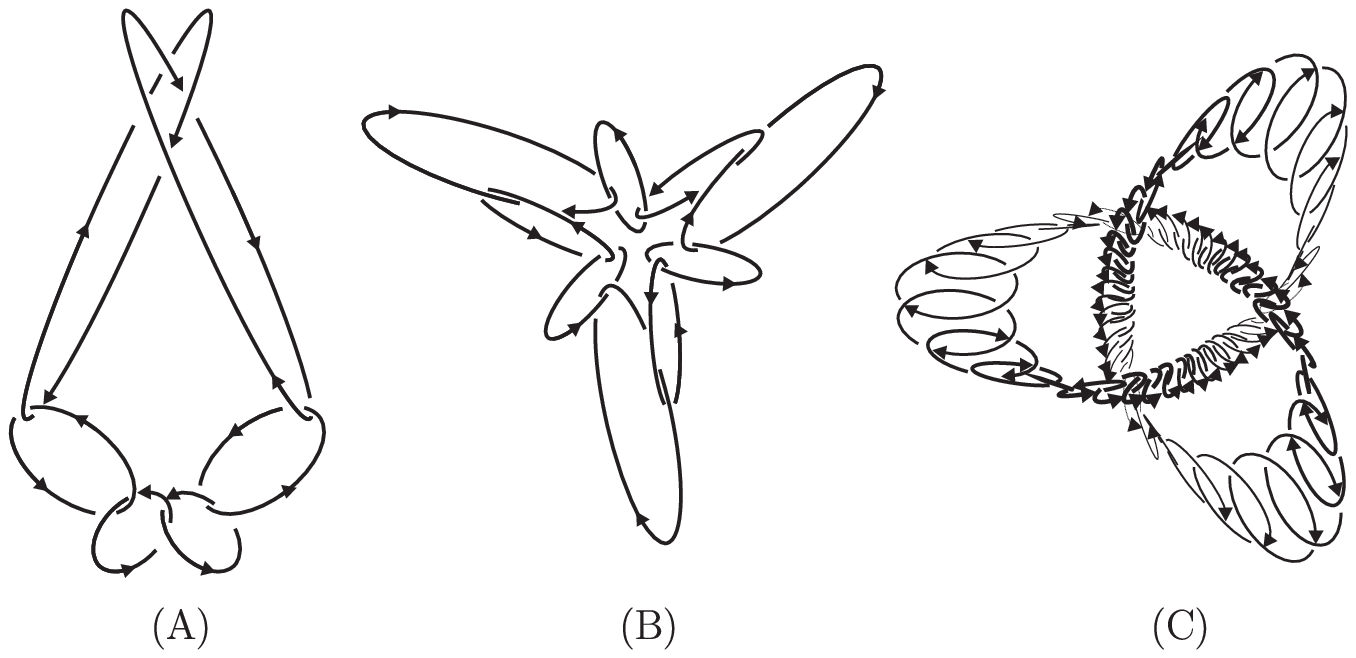}
\end{center}
\caption{(A)~$\Bd F(O,-1,6).$ (B)~$\Bd F(O,-3,9).$ 
(C)~$\Bd F(O\{2,-3\},-6,72)$.}
\label{more chains}
\end{figure}

Fig.~\ref{more chains} illustrates the application of
Theorem~\ref{Theorem A} to the knots in Fig.~\ref{legendrian 
knots}.

\section{Divides and fibered hyperbolic $\C$-line 
arrangements}\label{proof of Theorem B}

Following A'Campo \cite{A'Campo:divides97,A'Campo:divides00}, 
let a \emph{divide} be the image $\iota(J)$ of a non-empty
compact $1$-manifold-with-boundary $J$ by a smooth immersion 
$\iota: J\to D^2$ such that
\begin{enumerate}
\item $\iota$ is proper, that is, $\iota^{-1}(\Bd D^2)=\Bd J$;
\item $\iota$ is transverse to $\Bd D^2$;
\item $\iota(J)$ is connected;
\item $\iota(J)$ has only finitely many singular points, each of
which is a doublepoint with distinct tangent lines.
\end{enumerate}
A'Campo gives a simple, beautiful construction
(essentially a ``unoriented Gaussian resolution'' of $\iota$), 
phrased in terms of the (co)tangent bundle of $D^2$, which 
associates to any divide $\iota(J)\sub D^2$ a link $K(\iota(J))\sub S^3$. 
The components of $K(\iota(J))$ are in natural one-to-one correspondence 
with the components of $J$, and there is a natural orientation 
for $K(\iota(J))$ (up to simultaneous reversal of orientation on all 
components).  A'Campo then proves that $K(\iota(J))$ 
is a fibered link.

\begin{figure}
\begin{center}
\includegraphics[width=5in]{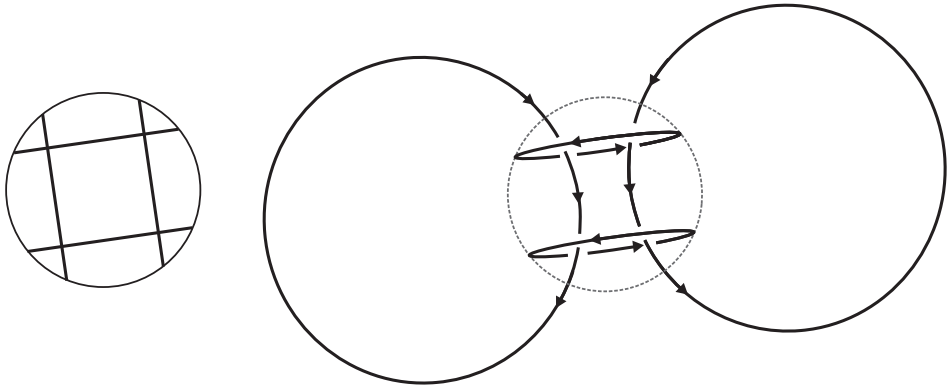}
\end{center}
\caption{A divide which is the closure 
of an essentially affine hyperbolic $\R$-line arrangement;
A'Campo's fibered link of this divide (which 
is also the link at infinity of the arrangement,
and isotopic to the link in Fig.~\ref{chain}), 
supplemented with the circle $\Bd\H\R\sub\Bd\H\C=S^3$.}
\label{square divide}
\end{figure}

Fig.~\ref{square divide} illustrates A'Campo's construction.  
(The image of the circle $\Bd\H\R\sub S^3$ by the stereographic 
projection in \ref{about Legendrian knots} lies in 
$\R\times\R\sub\C\times\R$, and thus projects 
via $\pr_1:\C\times\R\to\C$ onto a line segment;
in Fig.~\ref{square divide}, the diagram of the fibered 
link has been produced using the alternative projection 
$\mathrm{Re}\times\pr_2:\C\times\R\to\R\times\R$,
allowing it to be legibly supplemented with a diagram
of $\Bd\H\R$.)

\begin{thm}\label{Theorem B} The link at infinity of the
complexification $\C\Ls_\R\sub\Ps(\C)$ of an
$\R$-line arrangement $\Ls_\R\sub\Ps(\R)$ 
is fibered iff $\Ls_\R$ is connected.
\end{thm}

\begin{proof} If $\Ls_\R\sub\H\R$ is not connected, 
then $\Kinf{\C\Ls_\R}$ is a split link and so not
fibered.  If $\Ls_\R$ is connected, then (by Lemma~\ref{perturb})
there is a small perturbation $\Ls'_\R$ of $\Ls_\R$ such that
$\Kinf{\Ls'}$ and $\Kinf\Ls$ are ambient isotopic and $\Ls'_\R$
is in general position; of course $\Ls'_\R$ is also connected,
so its closure $\overline{\Ls'_\R}$ is a divide.  An inspection of
A'Campo's construction in \cite{A'Campo:divides97}, translating his 
symplectic approach into complex language, immediately
reveals that $K(\overline{\Ls'_\R})$ (with one of its two natural
orientations) and $\Kinf{\C\Ls'_\R}$ are identical, 
whence $\Kinf{\C\Ls_\R}$ is fibered.

The affine case reduces immediately to the hyperbolic case.  
\end{proof}

\section{Remarks}\label{remarks}

\subsection{Remarks on Section~\ref{hyperbolic line arrangements}}

(1)~As combinatorial structures, 
links at infinity of hyperbolic $\R$-line arrangements
are more or less identical to the ``chord diagrams'' 
used in the theory of Vassiliev invariants.
The significance of this coincidence, if any, 
is unclear.

(2)~The link at infinity $\Kinf\G \sub S^3$
of an affine $\C$-algebraic curve $\G\sub \C^2$, introduced in 
\cite{embeddings}, has been studied quite extensively 
(see, e.g., 
\cite{unfoldings,%
Neumann:curves_in_surfaces,%
Neumann:links_at_infinity,%
Ha:irregularity,%
Nemethi-Zaharia:Milnor_at_infinity,%
Neumann-Thanh,%
Cassou-Nogues-Ha:Lojasiewicz,%
Bodin}; see also 
\cite{Neumann-Norbury}).
There are strong topological restrictions on $\Kinf\G$,
and in particular $\Kinf\G$ is solvable
and ``nearly fibered'' (see, e.g., 
\cite{Neumann:links_at_infinity} and \cite{Neumann-Norbury}).
Typically, in $\C^2$, if $\G$ is an affine $\C$-analytic curve 
which is not $\C$-algebraic, then the isotopy type 
of $(1/r)(\G\cap S^3_r)\sub S^3$ does not stabilize
as $r\to\infty$, and so there is no apparent useful notion 
of the link at infinity of such a curve.
By contrast, in $\H\C$, where there is no obvious special class of 
$\C$-analytic curves which might usefully be called ``algebraic'', 
there is nonetheless 
a plentiful supply of $\C$-analytic curves (among them, 
curves equivalent by suitable diffeomorphisms $\H\C\to\C^2$ to 
arbitrary $\C$-algebraic curves in $\C^2$) with well-defined links 
at infinity.  
Such curves were studied from a topological viewpoint
in, e.g., \cite{someknot,totaltan,classical}, as
were the corresponding links (under such names as 
\emph{transverse $\C$-links}), which turn out to be
topologically much more diverse than in the affine case
(for instance, they are rarely solvable, and typically
they are very far from being fibered); 
however, the interpretation in terms of $\H\C$ has 
heretofore been completely overlooked.
I hope to come back soon to a new study
of general links at infinity in $\H\C$, in which the 
techniques of \cite{someknot} (etc.) will be supplemented with 
techniques of hyperbolic geometry.

\begin{figure}
\begin{center}
\includegraphics[width=5in]{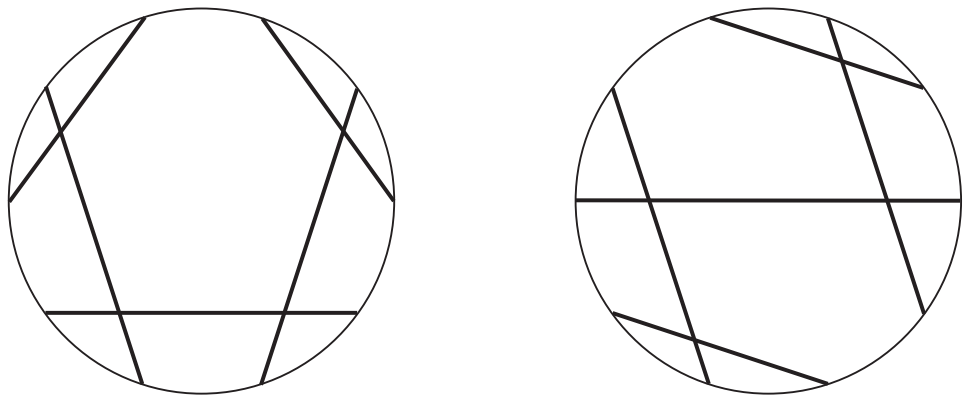}
\end{center}
\caption{There is no line in the hyperbolic $\R$-line arrangement 
$\Ls_\R$ (left) whose boundary separates the boundaries of two 
other lines in $\Ls_\R$.  
The boundary of the horizontal line in $\Ls'_\R$ (right) 
does separate the boundaries of two other lines in $\Ls'_\R$.  
Thus $\Kinf{\Ls_\R}$ and $\Kinf{\Ls'_\R}$ are not ambient isotopic,
though $\Kinf{\C\Ls_\R}=\Kinf{\C\Ls'_\R}$.}
\label{line tree}
\end{figure}

(3)~Hyperbolic $\R$-line arrangements with non-isotopic
links at infinity can easily have complexifications with
isotopic links at infinity.  For instance, if
$\Ls_\R=\bigcup_{i=1}^k L_i$ is a ``line 
tree''---that is, $\Ls_\R$ is in general position
and contractible (equivalently, the abstract $1$-complex 
$\G(\Ls_\R)$ dual to $\Ls_\R$, whose vertices are 
$L_1,\dots,L_k$ and whose edges are the nodes of 
$\Ls_\R$, is a tree)---then the ambient isotopy
type of $\Kinf{\C\Ls_\R}$ depends only on $\G(\Ls_\R)$ (it is
a ``$\G(\Ls_\R)$-connected sum'' of $k-1$ copies of $\Bd A(O,-1)$);
however, for $k>4$ there exist pairs $\Ls_\R$, $\Ls'_\R$ 
with $\Kinf{\Ls_\R}\ne \Kinf{\Ls'_\R}$ and 
$\G(\Ls_\R)=\G(\Ls'_\R)$.  
Fig.~\ref{line tree} gives an example.

\subsection{Remarks on Section~\ref{proof of Theorem A}}

(1) The word ``chain'' is already used in complex hyperbolic
geometry (see, e.g., \cite{Goldman}) to refer to certain 
important real curves in $S^3$ (and more generally $S^{2n-1}$,
considered as the boundary of complex hyperbolic $n$-space), 
which are not chains as defined here.  On the other hand,
chains as defined here (specifically, chains of type $O$) 
have been studied in the context of $3$-dimensional real 
hyperbolic geometry by Neumann and Reid \cite{Neumann-Reid}. 
Neither of these coincidences should cause confusion.

(2) It seems likely that the following
generalization of Theorem~\ref{Theorem A} can 
be established by similar techniques.

\begin{conj} For every quasipositive Seifert surface
$F\sub S^3$, there exists a system of pairwise disjoint proper
arcs $\a_1,\dots,\a_\mu\sub F$ such that, for all sufficiently
large $k_1,\dots,k_\mu\in\N$, if $F'$ is the result of plumbing 
$k_i$ copies of $A(O,-1)$ to $F$ along parallels of $\a_i$
for $i=1,\dots,\mu$, then $\Bd F'$ is isotopic to the link
at infinity of a hyperbolic $\C$-line arrangement.
\end{conj}
The second half of the proof of Theorem~\ref{Theorem A} would 
still apply, showing that for such a surface $F'$, $\Bd F'$ is 
fibered iff $\Bd F$ is fibered.

\subsection{Remark on Section~\ref{proof of Theorem B}}
A modified (and strengthened) 
version of Theorem~\ref{Theorem B} may be proved using 
other techniques (and other language), notably those of 
\cite{espaliers}.  
Let $\Ls_\R\sub\H\R$ be a connected hyperbolic $\R$-line 
arrangement.

\begin{thm}\label{Theorem C}
If $\Ls_\R\sub\H\R$ is essentially affine, 
then $\Kinf{{\CLs_\R}}$ is a closed positive braid.
In any case, there exists an espalier $\Ts$ such that 
$\Kinf{{\CLs_\R}}$ is a closed $\Ts$-positive braid.
\end{thm}
The proof will be given elsewhere.

\subsection{Further remarks}
As mentioned in Section~\ref{intro}, 
an \emph{open book} for an oriented link $K\sub S^3$
is a smooth map $f: S^3\to\C$ with $K=f^{-1}(0)$,
such that $0$ is a regular value of $f$ 
and $f/|f|:S^3\setminus f^{-1}(0)\to S^1$
is a fibration; and an oriented link $K$ is 
\emph{fibered} if there is an open book for $K$.  
A \emph{trivial unfolding}, 
as defined in \cite{unfoldings} (see also \cite{Kauffman-Neumann}),
is more or less the cone on an open book.
For example, if $H(k)$ is the fibered link comprised 
of $k\ge 1$ coherently oriented fibers of the standard 
positive Hopf fibration $S^3\to S^2$ (so $H(1)=O$ and
$H(2)=\Bd A(O,-1)$), then the open book 
$f_k:S^3\to\C:(z,w)\mapsto z^k+w^k$ for $H(k)$
yields the trivial unfolding 
$u_k:\{(z,w)\in D^4: |z^k+w^k|\le \e\} \to \C :(z,w)\mapsto z^k+w^k$
(for any sufficiently small $\e>0$).
The domain $D$ of a trivial unfolding $u$ is a smooth 
$4$-ball-with-corners, 
the corners of $D$ being the boundary of a tubular neighborhood 
$\Nb{{\Bd D}}{u^{-1}(0)}$,
and the restriction $u|\Bd D$ being an open book (modulo a
corner-smoothing identification of $\Bd D$ with $S^3$).

Generalizing the definition of trivial unfolding,
in \cite{unfoldings} (see also \cite{Neumann-Norbury}) 
an \emph{unfolding} is defined as a map $u:D\to\C$, where $D$ is a 
smooth $4$-ball-with-corners, $u|\Bd D$ is 
an open book (modulo corner-smoothing),
and $u|\Int D$ has only finitely many critical points 
$\mathbf x_i, i= 1, \dots, m$, near which 
it restricts to trivial unfoldings
(in appropriately adapted coordinates) $u_i:D_i\to\C$.
Write $K=\Unfold_{i=1}^{m} K_i$
to indicate that some such unfolding, 
with local fibered links $K_i$, exists;
this symbolic equation suppresses
the fact that $K_1,\dots, K_m$ by no means
determine $K$ (for $m>1$).  

Let $\Ls_\R$ be a connected hyperbolic $\R$-line arrangement with
$n$ nodes, of degrees $k_1,\dots, k_n$.  
Among the components of $\H\R\setminus\Ls_\R$,
let there be exactly $s$ which are hyperbolically bounded
(that is, have closure contained in $\H\R$).
Let $k_{n+1}=\dots=k_{n+s}=2$.

\begin{thm}\label{Theorem D} 
$\Kinf{\CLs_\R}=\Unfold_{i=1}^{n+s}H(k_i)$.
\end{thm}

For essentially affine $\Ls_\R$, this follows 
from \cite{unfoldings} (or \cite{Neumann-Norbury}).  
Various proofs can be given for arbitrary (connected) $\Ls_\R$:
one uses Theorem~\ref{Theorem C} in combination with techniques
of \cite{espaliers}; another proceeds (reasonably directly)
from A'Campo's results \cite{A'Campo:divides97,A'Campo:divides00}; 
a third, using the techniques (rather than just the language) of
knot-theoretical unfolding, and extending the
statement of the theorem considerably, will be given 
in \cite{unfolding_AC}.

\smallskip
\noindent
\textbf{Acknowledgements.} Thanks to the organizers of ``Arrangements in Boston''
for their invitation to participate in that fruitful conference.
Thanks to the referee for valuable comments on the exposition
of this paper.  Finally, thanks to Yuji Kida (author of UBASIC86)
and Keith Hertzer (author of Graphmatica), whose software made it
practicable to prepare Figs.~\ref{legendrian knots} and
\ref{more chains}.


\providecommand{\bysame}{\leavevmode\hbox to3em{\hrulefill}\thinspace}

\end{document}